\input amssym.def
\input amssym
\magnification=1200
\parindent0pt
\hsize=16 true cm
\baselineskip=13  pt plus .2pt
$ $

\def\Z{\Bbb Z}
\def\A{\Bbb A}

\def\R{\Bbb R}

\centerline {\bf On topological actions of finite groups on $S^3$}

\bigskip

\centerline {Bruno P. Zimmermann}

\bigskip

\centerline {Universit\`a degli Studi di Trieste} 
\centerline {Dipartimento di Matematica e Geoscienze} 
\centerline {34127 Trieste,  Italy}

\bigskip

{\bf Abstract.}  We consider orientation-preserving actions of a finite group $G$
on the 3-sphere $S^3$ (and also on Euclidean space $\Bbb R^3$). By the
geometrization of finite group actions on 3-manifolds, if such an action is smooth
then it is conjugate to an orthogonal action, and in particular $G$ is isomorphic
to a subgroup of the orthogonal group SO(4) (or of SO(3) in the case of $\R^3$). On
the other hand, there are topological actions with wildly embedded fixed point
sets; such actions are not conjugate to smooth actions but
one would still expect that the corresponding groups $G$ are isomorphic to
subgroups of the orthgonal groups SO(4) (or of SO(3), resp.). In the present paper,
we obtain  some results in this direction; we prove that the only finite, nonabelian
simple group with a topological action on $S^3$, or on any homology 3-sphere, is the
alternating or dodecahedral group $\A_5$ (the only finite, nonabelian simple
subgroup of SO(4)), and that every finite group with a topological,
orientation-preserving action on Euclidean space 
$\R^3$ is in fact isomorphic to a subgroup of SO(3).

\bigskip \bigskip

{\bf 1. Introduction}

\medskip

We consider topological actions of finite groups on the 3-sphere (i.e., by
homeomorphisms); all actions in the present paper will be faithful and
orientation-preserving.  By the geometrization of finite group actions on
3-manifolds after Thurston and Perelman, each finite group $G$ acting smoothly on
$S^3$ is conjugate to an orthogonal action; in particular, $G$ is isomorphic to a
subgroup of the orthogonal group SO(4). The last statement is no longer true for
smooth actions on arbitrary homology 3-spheres, and a
classification of such groups appears to be difficult (see [Z1], [MeZ1]). Also, for
each dimension $n \ge 4$ there  exist finite groups which admit a
topological, orientation-preserving action on
$S^n$ but are not isomorphic to a subgroup of SO(n+1) (see [Z3]; this remains
open for smooth actions).

\medskip

On the other hand, not much is known on the finite groups $G$ admitting a
topological action on $S^3$; some examples are known of actions with wildly
embedded fixed point sets, mainly for cyclic groups. Such actions are not conjugate
to smooth or orthogonal actions, but one would still expect that the groups $G$ are
isomorphic to subgroups of SO(4); in the present paper, we prove the following
results in this direction.

\bigskip

{\bf Theorem 1.}  {\sl  Let $G$ be a finite group of orientation-preserving
homeomorphism of $S^3$ or a homology 3-sphere such that each element of $G$ has
nonempty fixed point set; then $G$ is isomorphic to a subgroup of the orthogonal
group ${\rm SO}(3)$. In particular, this is the case if $G$ has a global fixed
point.}

\bigskip

By extending actions on $\R^3$ to its 1-point compactification $S^3$, this implies:

\bigskip

{\bf Corollary 1.}  {\sl   A finite group $G$ of orientation-preserving
homeomorphisms of Euclidean space $\R^3$ is isomorphic to a subgroup of ${\rm
SO}(3)$.}

\bigskip

The proof of Theorem 1 uses a purely algebraic
charaterization of the finite subgroups of the orthogonal group SO(3) by A. Miller
(Theorem 3 in section 3). In section 3 we shall present a short, direct proof of
Corollary 1 on the basis of Miller's result;
an independent proof is given in [KSu]. For smooth actions, Corollary 1 follows
from the geometrization of finite group actions on $\R^3$ (see [KSc]). In
particular, every smooth action of a finite group on
$\Bbb R^3$ has a global fixed point. This remains true for topological
actions of finite solvable groups on $\R^3$ (see [GMZ, proof of Theorem 1]); since
the only nonsolvable subgroup of ${\rm SO}(3)$ is the alternating group $\A_5$, this
raises the following:

\bigskip

{\bf Question.}  Is there a topological action of $\A_5$ on $\R^3$ without a global
fixed point? Equivalently, is there a topological action of $\A_5$ on $S^3$ with
exactly one global fixed point?

\bigskip

We note that the group $\A_5$ has an action on the Poincar\'e homology 3-sphere with
exactly one global fixed point; the complement of the fixed point is an acyclic
3-manifold with an action  of $\A_5$ without a global fixed point, so  the question
has a positive answer for homology 3-spheres and acyclic 3-manifolds, in general.

\bigskip

{\bf Corollary 2.} {\sl Let $M$ be a 3-manifold whose universal covering  is
$\R^3$ or $S^3$. Let $G$ be a finite group of orientation-preserving homeomorphisms
of $M$ with a global fixed point; then $G$ is isomorphic to a subgroup of ${\rm
SO}(3)$.}

\bigskip

This  follows from Corollary 1 and Theorem 1 by lifting $G$ to
an isomorphic group of homeomorphisms of the universal covering with a
global fixed point. Note that, by considering invariant regular neighbourhoods,
this is true for smooth actions for the case of arbitrary 3-manifolds, and also
for topological actions it should remain true for arbitrary 3-manifolds.

\bigskip

{\bf Theorem 2.} {\sl   A finite, nonabelian simple group of homeomorphisms of
$S^3$ or a homology 3-sphere is isomorphic to the alternating group  $\A_5$.}

\bigskip

For smooth actions, this is proved in [Z2] and [MeZ2];  using Proposition 1 in
section 2, we reduce the proof of Theorem 2 to the proof in [MeZ2], based on
the Gorenstein-Walter classification of the finite simple groups with dihedral Sylow
2-subgroups. Actions of finite simple groups on spheres and homology spheres in
higher dimensions are considered in [GZ], and again these results remain true for
topological actions.

\medskip

In section 2 we prove a technical key result which allows to
generalize various  known results about finite group actions on homology 3-spheres
from smooth actions to topological actions.

\bigskip

{\bf 2. A preliminary result}

\medskip

Our technical key result is the following:

\bigskip

{\bf Proposition 1.}  {\sl    Let $G$ be a finite group with  an
orientation-preserving, topological action on $S^3$ (or on a homology 3-sphere).
Suppose that the fixed point set of an element $g \in G$ is a (possibly wildly
embedded) circle $K \cong S^1$. Then the normalizer $N = N_G(g)$ of $g$ in $G$ is
isomorphic to a subgroup of a semidirect product  $(\Z_a \times \Z_b) \rtimes \Z_2$
where $\Z_2$ acts dihedrally on the abelian group $A = \Z_a \times \Z_b$.}

\bigskip

For smooth actions this follows easily, for the case of arbitrary 3-manifolds, from
the existence of invariant regular neighbourhoods (considering rotations
of minimal angle around a smoothly embedded circle, see [MeZ2, Lemma 1]), and also
for topological actions it should remain true for arbitrary 3-manifolds. As in
the case of smooth actions, Proposition 1 has various applications to the structure
of finite groups admitting a topological action on $S^3$ (it is a basic tool for
the partial characterization of the finite groups which admit a smooth action on a
homology 3-sphere in [MeZ1], [Z1]).

\bigskip

{\it Proof of Proposition 1.}  Each element of the normalizer $N$ maps the circle
$K$ and its complement
$M = S^3 - K$ to itself. By Alexander-Lefschetz duality, $S^3 - K$ is also a
homology (and cohomology) circle (that is, has the homology or cohomology of the
circle).  Let
$B$ be the subgroup of $N$ which fixes $K$ pointwise;  then, by Smith fixed point
theory, $B$ acts freely on the cohomology circle $S^3 - K$, therefore  $B$ has
periodic cohomology of period two and is a cyclic group $B = \Z_n$  (see [Br]). We
call the elements of $B = \Z_n$ {\it rotations around} $K$. The factor group $C =
N/B$ acts faithfully on the circle $K$ and hence is a cyclic or a dihedral group;
its nontrivial elements are either {\it reflections} of $K$ (i.e., fixing exactly
two points of $K$), or {\it rotations along}
$K$, i.e. without fixed points on $K$.

\medskip

The normalizer $N$ acts on $S^3 - K$, with the cyclic normal subgroup $B = \Z_n $
acting freely, and we consider the cyclic regular covering $M = S^3-K  \to \bar M =
(S^3-K)/B$ and the cohomology spectral sequence associated to this covering (see
[McL]):  
$$E_2^{i,j}  =  H^i(\Z_n; H^j(M))   \;\;\;  \Rightarrow  \;\;\; H^{i+j}(\bar M),$$ 
converging to the graded group associated to a filtration of
$H^*(\bar M)$ (integer coefficients). Since $M$ is a cohomology circle, the
spectral sequence is concentrated in the rows $j=0$ and $j=1$, so the only possibly
nontrivial differentials, of bidegree $(2,-1)$, are 
$$d_2^{i,1}:  E_2^{i,1} = H^i(\Z_n; H^1(M))  \;\;  \to   \;\;   E_2^{i+2,0} =  
H^{i+2}(\Z_n; H^0(M)),$$

where $B = \Z_n$ acts  by the identity on $H^0(M) \cong \Z$ and, by duality, also
on $H^1(M) \cong \Z$ (since $G$ acts orientation-preservingly on $S^3$ and also on
$K$, it acts by the identity on the first homology and cohomology $\Z$ of $K$, and
hence by duality also on  the first homology and cohomology $\Z$ of $M = S^3 - K$). 
Summarizing, we have
$$d_2^{i,1}: H^i(\Z_n; \Z)  \;\;  \to   \;\;  H^{i+2}(\Z_n; \Z),$$

where $H^i(\Z_n;\Z)$ is isomorphic to $\Z$ if
$i=0$, to $\Z_n$ if $i>0$ is even, and trivial if $i$ is odd. In particular, $\Z_n$
has periodic cohomology of period two (for $i>0$), and the duality isomorphism is
given by the cup-product with an element $u \in H^2(\Z_n;\Z) \cong \Z_n$ (see [Br]).

\medskip

Since $\bar M$ is also a 3-manifold, its cohomology is trivial in dimensions larger
than three; hence, passing to the limit of the spectral sequence, the
differentials  $d_2^{i,1}$ have to be isomorphisms for larger $i$. But then also
the differential $d_2^{1,1}$ has to be an isomorphism, by dimension shifting with
the cup-product with $u \in H^2(\Z_n;\Z)$ and the multiplicative structure of the
spectral sequence (see [S, section 9.4]).   Replacing $H^0(\Z_n;\Z) \cong \Z$ by
its quotient, the Tate cohomology group $\hat H^0(\Z_n; \Z) \cong \Z_n = \Z/n\Z$,
also $d_2^{0,1}$ becomes an isomorphism, and hence $d_2^{0,1}: H^0(\Z_n; \Z) \cong
\Z  \to  H^2(\Z_n;
\Z) \cong \Z_n$ has to be surjective.  Passing to the limit of the spectral
sequence, this implies easily that $\bar M$ is a cohomology circle, and hence  also
a homology circle; in particular, $H_1(\bar M,\Z) \cong Z$.

\medskip

We consider the group extension 
$1 \to \pi_1(M) \hookrightarrow \pi_1(\bar M) \to \Z_n \to 1$ associated to the
regular covering $M \to \bar M$, with covering group $\Z_n$. Abelianizing
$\pi_1(M)$, we get an extension 
$$0 \to \pi_1(M)_{\rm ab} = \pi_1(M)/[\pi_1(M), \pi_1(M)] 
\hookrightarrow \pi_1(\bar M)/[\pi_1(M), \pi_1(M)]  \to \Z_n \to 0,$$ with 
$\pi_1(M)_{\rm ab} \cong H_1(M) \cong \Z$.  Since, as before, the covering group
$\Z_n$ acts trivially on $H_1(M) \cong \Z$, also   
$\pi_1(\bar M)/[\pi_1(M), \pi_1(M)]$ is an abelian group and hence isomorphic to the
abelianized group  $\pi_1(\bar M)_{\rm ab} \cong H_1(\bar M) \cong \Z$, so we have
an exact sequence   
$$0 \to H_1(M) \cong Z  \to H_1(\bar M) \cong \Z \to  \Z_n \to 0.$$ Let $g \in N$
be a rotation along $K$; since $g$ acts orientation-preservingly on
$K$,  it acts by the identity on the homology and cohomology of $K$, and by duality
by the identity also on $H_1(M) \cong \Z$. Then also the homeomorphism $\bar g$
induced by $g$ on
$\bar M$ acts by the identity on $H_1(\bar M) \cong \Z$, and hence $g$ acts by the
identity on the covering group $\Z_n$.  It follows that the subgroup $A$ of
$N$ of rotations around and along $K$ is abelian, of rank one or two. 

\medskip

If $g \in N$ is a reflection of $K$ instead, then $g$ fixes exactly two points in
$K$ and, by Smith theory, a circle in $S^3$  which implies $g^2 = 1$ (since a
nontrivial element cannot fix two different circles).  Also, since $g$ acts
dihedrally on the homology and cohomology of $K$, by duality it acts dihedrally also
on $H_1(M) \cong \Z$, and hence dihedrally also on $\Z_n$.

\medskip

Concluding, the normalizer $N$ has the structure given in Proposition 1.

\bigskip

{\bf 3.  Proofs of Theorems 1 and 2}

\medskip

The proof of Theorem 1 uses the following purely algebraic
characterization  of the finite subgroups of the orthogonal group
${\rm SO}(3)$.

\bigskip

{\bf Theorem 3.} ([Mi])  {\sl A finite group $G$ is isomorphic to a subgroup of 
${\rm SO}(3)$ if and only if the normalizer of each nontrivial element $g \in G$ is
a cyclic or dihedral group. This is in turn equivalent to the following condition: 
Distinct maximal cyclic subgroups of $G$ intersect trivially, and   each maximal
cyclic subgroup of $G$ has index one or two in its normalizer in $G$.}

\bigskip

{\it Proof of Theorem 1.}  By Smith fixed point theory, the fixed point set of a
nontrivial element $g \in G$ is a circle $K$. By Proposition 1, the normalizer
$N = N_G(g)$ is isomorphic to a subgroup of a semidirect product 
$(\Z_a \times \Z_b) \rtimes \Z_2$ where $\Z_2$ acts dihedrally on the abelian group
$A = \Z_a \times \Z_b$. As in the proof of Proposition 1, the group
$A$ consists of rotations around and along $K$; it has  a cyclic subgroup $B$ of
rotations around $K$ (i.e., fixing $K$ pointwise), with cyclic factor group $A/B$.
We will show that $A$ is a cyclic group and apply Theorem 3.

\medskip

Suppose that $A$ has a subgroup $\Z_p \times \Z_p$, for a prime $p$. We apply the
Borel formula (which holds in a purely topological setting) to the subgroup $\Z_p
\times \Z_p$ ([Bo]; see also [MeZ2] for such an application). If $p > 2$ then by the
Borel formula, $\Z_p \times \Z_p$ has exactly two cyclic subgroups $\Z_p$ with
nonempty fixed point set (two different circles), so some subgroup $\Z_p$ acts
freely on $S^3$ contrary to the hypothesis of Theorem 1 i). If $p=2$ then either
there are again two involutions in 
$\Z_2 \times \Z_2$ with nonempty fixed point set and one free involution, or all
three involutions have nonempty fixed point set (three circles intersecting in two
points). In particular, in the second case some involution in $\Z_2 \times \Z_2$
acts as a reflection on $K$ which, however, is not the case in the group $A$ of
rotations  of  $K$.

\medskip

So $A$ has no subgroups $\Z_p \times \Z_p$ and hence is cyclic, hence the normalizer
of $g$ in $G$ is cyclic or dihedral; by Theorem 3,  $G$ is isomorphic  to a
subgroup of ${\rm SO}(3)$, concluding the proof of Theorem 1.

\bigskip

{\it Proof of Corollary 1.}  We give a short direct proof of Corollary 1 which  uses
only a small part of the preceding results.

\medskip

The action of $G$ on $\R^3$ extends to an action on its 1-point
compactification $S^3$ with a global fixed point. We shall verify the 
conditions of the second part of Theorem 3 for this action of $G$ on $S^3$. 

\medskip

Let $g$ be an element of $G$ generating a maximal cyclic subgroup $M$ of $G$. By
Smith fixed point theory, the fixed point set of $g$ and $M$ is a single circle
$K$.  The normalizer  $N$ of $g$ in $G$ maps $K$ to itself, hence every nontrivial
element of $N$ fixes
$K$ pointwise  or acts as a reflection on $K$, and the subgroup $B$ of $N$ fixing
$K$ pointwiss has index one or two in $N$. By the first paragraph of the proof of
Proposition 1,  $B$ is cyclic; since $M$ is maximally cyclic, $B = M$, and hence $M$
has index one or two in its normalizer.

\medskip

Moreover, if  two maximal cyclic subgroups of  $G$ have nontrivial intersection
then they have the same circle $K$ as fixed point set, hence generate a cyclic
subgroup of $G$ and are equal. Theorem 3 now implies that $G$ is isomorphic to a
subgroup of SO(3).

\medskip

This concludes the proof of Corollary 1.

\bigskip

Finally, the {\it Proof of Theorem 2} is analogous to the proof of the main Theorem
in [MeZ2] (which uses the Gorenstein-Walter classification of the
finite simple groups with dihedral Sylow 2-subgroups). The proof of [MeZ2, Theorem]
depends on two lemmas;   Lemma 1 in [MeZ2] is the version for
smooth actions of Proposition 1 of the present paper, and Lemma 2 is
again a consequence of the Borel formula, so both Lemma 1 and Lemma 2 of [MeZ2] hold
for purely topological actions. The proof of Theorem 2 is now completely analogous
to the proof of [MeZ2, Theorem]: by Gorenstein-Walter one reduces first to the
linear fractional groups ${\rm PSL}_2(q)$, for an odd prime power $q$, or the
alternating group $\A_7$ (these are the finite, nonabelian simple groups with
dihedral Sylow 2-subgroups), and finally to $\A_5 \cong {\rm PSL}_2(5)$ by
Proposition 1 and the purely topological argument in [Z2] (using also [MeZ1, section
6] to exclude the small groups ${\rm PSL}_2(5^2)$ and ${\rm PSL}_2(3^2) \cong
\A_6$).

\bigskip  \bigskip

\centerline {\bf References}

\bigskip

\item {[Bo]} A. Borel, {\it Seminar on Transformation Groups.} Annals of Math.
Studies 46, Princeton University Press 1960

\item {[Br]} K.S. Brown, {\it Cohomology of Groups.}  Graduate Texts in Mathematics
87, Springer 1982

\item {[GMZ]} A. Guazzi, M. Mecchia, B. Zimmermann, {\it On finite groups acting on
acyclic low-dimensional manifolds,}  Fund. Math. 215,  203-217  (2011)

\item {[GZ]} A. Guazzi, B. Zimmermann, {\it On finite simple groups acting on
homology spheres.}  Mo-natsh. Math. 169,  371-381 (2013)

\item {[KSc]} S. Kwasik, R. Schultz, {\it Icosahedral group actions on $S^3$.}
Invent. Math. 108, 385-402  (1992)

\item {[KSu]} S. Kwasik, F. Sun, {\it Topological symmetries of $\R^3$.}
arXiv:160205295v2

\item {[McL]} S. MacLane, {\it Homology.} Springer-Verlag,  Berlin 1963

\item {[MeZ1]} M. Mecchia, B. Zimmermann, {\it On finite groups acting on
$\Z_2$-homology 3-spheres.} Math. Z. 248, 675-693 (2004)

\item {[MeZ2]} M. Mecchia, B. Zimmermann, {\it On finite simple groups acting on
integer and mod 2 homology 3-spheres.}  J. Algebra 298, 460-467  (2006)

\item {[Mi]} A. Miller, {\it A group theoretic characterization of the
2-dimensional spherical groups.}  Canad. Math. Bull. 32,  459-466 (1989)

\item {[S]} E.H. Spanier, {\it Algebraic Topology.}  McGraw-Hill, New York 1966

\item {[Z1]} B. Zimmermann, {\it On the classification of finite groups acting on
homology 3-spheres.} Pacific J. Math. 217, 387-395  (2004)

\item {[Z2]} B. Zimmermann, {\it On finite simple groups acting on homology
3-spheres.} Top. Appl. 125, 199-202 (2002)

\item {[Z3]} B. Zimmermann, {\it On topological actions of finite, non-standard
groups on spheres.} arXiv:1602.04599 (to appear in Monatsh. Math.)

\item {[Z4]} B. Zimmermann, {\it On finite groups acting on spheres and finite
subgroups of orthogonal groups.}  Sib. Electron. Math. Rep. 9, 1-12 (2012)
(http://semr.math.nsc.ru)

\bye